\documentclass[12pt, a4paper]{article}
\usepackage{amsmath}
\usepackage{amsfonts}
\usepackage{amsthm}
\usepackage[latin1]{inputenc}
\newtheorem{theorem}{Theorem}[section]

\newtheorem{lemma}[theorem]{Lemma}

\newcommand{\supp}{\operatornamewithlimits{supp}}
\numberwithin{equation}{section}
\newcommand{\eref}[1]{{\rm (\ref{#1})}}
\newcommand{\lref}[1]{{\rm Lemma~\ref{#1}}}
\newcommand{\tref}[1]{{\rm Theorem~\ref{#1}}}

\begin{document}
\title{A counterexample to a weak-type estimate for potential spaces and tangential approach regions}
\author{\textbf{Javier Soria\footnote{Research partially supported by Grants BFM2001-3395 and  
2001SGR00069.}}\\ 
Dept. Appl. Math. and Analysis\\
University of Barcelona\\
E-08071 Barcelona, SPAIN\\
E-mail: soria@mat.ub.es\vspace{0.2cm}\\
\textbf{Olof Svensson}\\
Dept. of Science and Technology\\
Campus Norrköping, 
Linköping University\\
SE-601 74 Norrköping, SWEDEN\\
E-mail: olosv@itn.liu.se
}
\maketitle
\begin{abstract}\noindent
  We show that for every nontrivial potential space $L^{1}_{K}(\mathbb{R}^{n})$, there exists an approach
  region for which the associated maximal function is of weak-type, but the
  boundedness for the completed region is false, which is in contrast
  with the nontangential case.
\end{abstract}\vspace{0.5cm}
\vspace{0.5cm}

\section{Introduction}
In \cite{NS84} it was proved that Fatou's theorem holds on regions
$\Omega$, larger than cones (but still nontangential), by means of the
boundedness of the associated maximal function $M_{\Omega}$. One of
the key points in that proof is that one could replace the given
region, by a larger region $\widehat{\Omega}$ obtained by adding a
cone at any point of $\Omega$, and then prove that the boundedness of
the two maximal functions $M_{\Omega}$ and $M_{\widehat{\Omega}}$ 
are equivalent. This seems geometrically very natural, since the difference, at any point, between
$\widehat{\Omega}$ and $\Omega$, is just the canonical approach region (i.e., a cone).

In \cite{NRS82} Fatou's theorem was extended to some tangential approach
regions, when the functions were assumed to have some a priori
smoothness (they belonged to a potential space). This result was later
on generalized in \cite{RS} to characterize all the approach regions
(under a completion hypothesis similar to the one in \cite{NS84}) for which the convergence holds for
the potential spaces. 

The main result of this paper is to show that, contrary to the case of
\cite{NS84}, the assumptions on the region assumed in \cite{RS}, which is
natural as we mentioned before, from the point of view of convergence, turns out to give
different boundedness results for the corresponding maximal
operators. In order to clarify this statement, let us introduce some
notations: 
\medskip

Let $P_{t}(x)$ be the Poisson kernel in $\mathbb{R}^{n+1}_{+}$.
Given a set $\Omega \subset\mathbb{R}^{n+1}_{+}$, we define the 
 maximal function: 
\[
M_{\Omega}f(x)=\sup_{(y,t)\in\Omega_{x}}
|P_{t}*f(y)|,
\]
where $\Omega_{x}=x+\Omega$.

If $r:\mathbb{R}^{+}\to\mathbb{R}^{+}$ is
an increasing function, then we define the \lq\lq cone" for the function $r$ as:
\[
\Gamma_{r}(x,t)=\{(y,s):|x-y|\leq r(s)-r(t)\}.
\]
If $r(t)=t$, then $\Gamma_{t}=\Gamma$ is a nontangential cone.

We say that $\Omega$ satisfies the $r$-condition if
$\Gamma_{r}(x,t)\subset\Omega$ for all $(x,t)\in\Omega$.
 For example, in the case of nontangential
approach, $r(t)=t$ and the $r$-condition is the cone condition of
\cite{NS84}. The function $r$ is determined, in each case, from the
potential space under consideration. 

The potential space considered here is $L^{1}_{K}(\mathbb{R}^{n})$:
\[
L^{1}_{K}(\mathbb{R}^{n})=\{f:f=K*F,\ F\in L^{1}(\mathbb{R}^{n})\}.
\]
The  kernel $K$ is positive and integrable, but unbounded ($K(0)=\infty$), 
nonnegative and radial  (if
$|x|=|y|$, then $K(x)=K(y)$), and decreasing  (if $|x|\leq |y|$, then
$K(x)\geq K(y)$).   
We consider the following  norm on the potential space $L^{1}_{K}(\mathbb{R}^{n})$:
\[
\|f\|_{L^{1}_{K}(\mathbb{R}^{n})}=\inf_{f=K*F}\|F\|_{L^{1}(\mathbb{R}^{n})}.
\]

For the space $L^{1}_{K}(\mathbb{R}^{n})$, we have that if  
$r_{K}(t)=\|P_{t}*K\|_{\infty}^{-1/n}$,  
then the region
$\Gamma_{K}=\Gamma_{r_{K}}$  is tangential, under the above 
assumptions on the kernel $K$  (see \cite{NRS82}). This can be expressed as
\begin{equation}\label{tang}
\lim_{t\to 0}\frac{r_{K}(t)}{t}=\infty.
\end{equation}

In case of the Bessel potential spaces 
$L^{1}_{\alpha}({\mathbb{R}^{n}})=\{F*G_{\alpha}:F\in L^{1}(\mathbb{R}^{n})\}$
 (where $G_{\alpha}$ is the Bessel potential), then 
$r_{G_{\alpha}}(t)=t^{1-\alpha /n}$.

As a consequence of Theorem 2.6 in \cite{RS}, we know that if
$\Omega$ satisfies the $r_{K}$-condition, then 
$M_{\Omega}:L^{1}_{K}(\mathbb{R}^{n})\to L^{1,\infty}(\mathbb{R}^{n})$ 
if and only if $|\Omega (t)|\leq C(r_{K}(t))^{n}$, 
for all $t>0$, where 
$\Omega (t)=\{x:(x,t)\in\Omega\}$.

Given an approach region $\Omega$, we can always define the smallest
region containing $\Omega$, satisfying the $r_{K}$-condition as
follows:
\[
\widehat{\Omega}_{K}=\{(y,t)\in\mathbb{R}^{n+1}_{+}:
\exists (x,s)\in\Omega,|x-y|\leq r_{K}(t)-r_{K}(s)\}.
\]
Then it is easy to show that $\widehat{\Omega}_{K}$ satisfies the
$r_{K}$-condition, and $\Omega\subset\widehat{\Omega}_{K}$.

In the nontangential case it was proved in \cite{NS84} that
$M_{\Omega}:L^{1}(\mathbb{R}^{n})\to L^{1,\infty}(\mathbb{R}^{n})$ if
and only if  
$M_{\widehat{\Omega}}:L^{1}(\mathbb{R}^{n})\to L^{1,\infty}(\mathbb{R}^{n})$.
However, we will show in \tref{pr} that under the above conditions on $K$, and hence (\ref{tang}) holds,
then this equivalence fails in general. This is somehow surprising, since 
$M_{\Gamma_{K}}:L^{1}_{K}
(\mathbb{R}^{n})\to L^{1,\infty}(\mathbb{R}^{n})$
(see \cite{NRS82}).  Therefore, even though the boundary convergence holds within both $\Omega$ and the
\lq\lq cone" $\Gamma_K$, it fails for the completed region $\widehat{\Omega}_K$.

\section{Main theorem}

We now prove our main result, namely  that the characterization in \cite{NS84} does
not hold for tangential regions:  
a maximal operator $M_{\Omega}$ can be of weak type (1,1) while the
maximal operator for the completed region, $M_{\widehat{\Omega}_{K}}$ fails 
to be of weak  type (1,1).

\begin{theorem}\label{pr}
  For each of the potential spaces  $L_{K}^{1}(\mathbb{R}^{n})$,
  there   exists a region   
  $\Omega$ with the following properties:
  \item{(i)} $\Omega$ satisfies the cone condition.
  \item {(ii)} $|\Omega(t)|\leq C(r_{K}(t))^{n}$.
  \item {(iii)} $|\{M_{\Omega} f>\lambda\}|\leq
    C\displaystyle\frac{\|f\|_{L^{1}_{K}}}{\lambda}$.
  \item {(iv)}  $\displaystyle\frac{|\widehat{\Omega}_{K}(t)|}{(r_{K}(t))^{n}}$ is
    unbounded. 
  \item  {(v)} $M_{\widehat{\Omega}_{K}}$ is not of weak    type (1,1).
\end{theorem}
\noindent
The proof  uses the following  lemma from \cite{Sj83}.

\begin{lemma}\label{kbit}
  Assume the operators $T_{k}$, $k=1,2, \dots ,$ are defined
  in  ${\mathbb{R}}^{n}$ by
  \begin{equation}\label{kb4}
    T_{k}f(x)=\sup _{s\in I_{k}}(K_{s}*|f|(x)),
  \end{equation}
  where the $K_{s}$ are integrable and non-negative in
  ${\mathbb{R}}^{n}$, and the index sets $I_{k}$ are such that $T_{k}f$
  are measurable for any measurable $f$. Let for each $i=1,\dots ,n$ a
  sequence $(\gamma _{ki})_{k=1}^{\infty }$ be given with $\gamma
  _{ki}\geq \gamma _{k+1,i}>0$ and assume the $T_{k}$ are uniformly of
  weak type (1,1), and
  \begin{equation*}
    \supp K_{s}\subset \{x=(x_{1},\dots ,x_{n})\in
    {\mathbb{R}}^{n};|x_{i}|\leq \gamma _{ki} ,\ i=1,\dots ,n\},\ s\in I_{k},
  \end{equation*}
  and
  \begin{equation*}
    \int K_{s}^{*}\leq C_{0},\ \ s\in \bigcup _{k} I_{k},
  \end{equation*}
  where for $s\in I_{k}$
  \begin{equation}\label{kb3}
    K_{s}^{*}(x)=\sup \{K_{s}(x+y);|y_{i}|\leq \gamma _{k+N,i},\
    i=1,\dots ,n\}
  \end{equation}
  for some fixed natural number $N$. Then the operator 
  \[Tf(x)=\sup _{k}T_{k}f(x),\]
  is of weak type (1,1).
\end{lemma}

\begin{proof}[Proof of Theorem \ref{pr}] For simplicity, we  will usually drop the subscript $K$, and we
will write
  $r(t)=r_{K}(t)$, although for the regions $\Gamma_{K}$ we will keep it. Also, we only consider the
case $n=1$ (higher dimensions require minor modifications). 

  We start with the construction of the region $\Omega$: for this we
  choose a set of points $\omega$ from which we obtain the 
  region $\Omega$ by completing $\omega$ with nontangential cones.
  
  To construct $\omega$, we define a curve $\gamma (t)$, 
  \[
  \gamma (t)=N(t)r(t),
  \]
  where $N(t)$ is an integer valued function  that tends to infinity
  as $t\to 0$. There are some restrictions on how fast $N(t)$ may
  increase, which will be explained below.    

  The first condition on $N(t)$ is that $\gamma (t)=N(t)r(t)$
  should tend to zero as $t\to 0$. 
  We also need a sequence $t_{k}$ decreasing to $0$ fast
  enough, the precise meaning of this is described later on.
  
  Consider the tangential curve $(\gamma (t),t)$, for small $t>0$.
  If we compare the curve $\gamma (t)$ to the curve $r(t)$, which
  defines the region $\Gamma_{K}$, 
  we get
  \begin{equation}
    \frac{\gamma (t)}{r(t)}=N(t)\to\infty\ \mathrm{as}\ t\to
    0.
  \end{equation}
  This shows that $\gamma$ is well outside $\Gamma_{K}$. 
  
  Now choose a starting level, $t_{1}$, this will need to be small,
  exactly how small, will be made clear below.  Let  $x_{1}^1=\gamma
  (t_{1})$. 
  The first $N(t_{1})$ points in the set $\omega$ are
  \[
  \{(x_{i}^1,t_{1}):\ 
  x_{i}^1=x_{1}^1-ir(t_{1}),\ 0\leq i\leq N(t_{1})-1\}.
  \]
  These points  $(x_{i}^1,t_{1})$ have to be well 
  outside the nontangential cone $\Gamma$. More precisely:
  \[
  \gamma(t_{1})-(N(t_{1})-1)r(t_{1})=r(t_{1})>3t_{1}.
  \]
  If $t_{1}$ is small enough this is true, due to the tangentiality of 
  $\Gamma_{K}$ (see (\ref{tang})), and we choose $t_{1}$ to be any such number.  
  
  We proceed inductively, assuming that we have chosen $t_{k-1}$, and
  added the $N(t_{k-1})$ points at this level to $\omega$. 
  
  Now choose any $t_{k}<t_{k-1}$ satisfying
  \begin{equation}
  \label{con}
  \gamma (t_{k})+(t_{k-1}-t_{k})<2t_{k-1}
  \end{equation}
  which is to say that after adding the nontangential cone to
  $(\gamma(t_{k}),t_{k})$ the region thus obtained is contained in the
  nontangential cone $\Gamma_{2t}$, at height $t_{k-1}$. It is  
  obvious that this   cone does not intersect  the previously chosen
  points in $\omega$.    
  Now add the following $N(t_{k})$ points to the set $\omega$:
  \[
  \{(x_{i}^k,t_{k}):\ x_{i}^k=\gamma(t_{k})-ir(t_{k}),\ 0\leq
  i\leq N(t_{k})-1\}.
  \]
  This finishes the construction on the level $t_{k}$. If we  continue 
  this way, the set of points $\omega$ is obtained.

  We end up with a set of points $\omega$ which are arbitrarily
  close to the boundary, whose number   at
  height $t_{k}$   increases  to infinity as $t_{k}\to 0$.
  The region $\Omega$ is then obtained   by completing
  $\omega$ with the nontangential cone $\Gamma$.
  
  We still have to determine the function $N(t)$: we impose also another condition
     to make 
  sure that the cross-sections of $\Omega$ will satisfy the right estimate: 
  \begin{equation}\label{size}
    |\Omega (t)|\leq Cr(t).
  \end{equation}
  We start at any level $t_{k}$ and move upwards to $t_{k-1}$. At the
  level $t_{k}$ the region $\Omega$ consists of one part that is
  contained in a fixed nontangential cone with vertex at the origin, 
  this part comes from the
  lower levels (see \eref{con}), and here the size estimate \eref{size} is obvious.

  The other part consists of $N(t_{k})$
  points, as we move upwards, first  each interval will have  a size at
  height $t$ 
  which is bounded from above by $t$, this is the case if $t$ is below 
  $t_{k}+\frac{1}{2}r(t_{k})$, the size of the union of these
  intervals is then bounded from above by $tN(t_{k})$. For larger $t$ the
  intervals will 
  have met and the size estimate follows, if it holds while they are
  disjoint. 
  If we impose on $N(t)$ that
  \begin{equation}
    \label{eq:cond}
    tN(t)< r(t)
  \end{equation}
  then the size of the disjoint intervals will have the correct  upper
  bound. Since the region  $\Gamma_{K}$ is tangential,   
  this can be achieved, while $N(t)$ tends to infinity as $t\to 0$.

  If we instead complete the region $\Omega$ with the tangential
  region associated with the   potential space $L_{K}^{1}$, that is
   ${\Gamma}_{K}$, then  
   $|\widehat{\Omega}_{K}(t)|/r(t)$ will not be
   bounded, since otherwise we could find a constant $C$ such that:

   \begin{equation}
     \label{compsize}
     |\widehat{\Omega}_{K}(t)|\leq Cr(t).
   \end{equation}
   At levels $t\in [t_{k},2t_{k}]$ we will have $N(t_{k})$
   intervals that are almost disjoint (observe that by (\ref{con}), $2t_k<t_{k-1}$), and the measure of
the union of
   these intervals can be estimated from below by a constant times 
   $N(t_{k})r(t_{k})$ at height $2t_{k}$. We thus have 
   \[
   |\widehat{\Omega}_{K}(2t_{k})|\geq CN(t_{k})r(t_{k}).
   \]
   Letting $t=2t_{k}$ in \eref{compsize} we see that in order for the
   above estimates  to be compatible, we must have
   \[
   N(t_{k})r(t_{k})\leq Cr(2t_{k})
   \]
   and this is only possible if $N(t)$ is bounded, since
   $r(t_{k})\sim r(2t_{k})$ (this 
   follows from the related relation for the Poisson kernel,
   $P_{t}(x)\sim P_{2t}(x)$). Therefore, $\widehat{\Omega}_{K}$ cannot
   satisfy the necessary condition \eref{size}, and hence,  
   $M_{\widehat{\Omega}_{K}}$ cannot be of weak type (1,1) (by Theorem 2.6 in\cite{RS}).
   
   Now that the region $\Omega$ is defined (and we have dealt with (i), (ii), (iv), and (v) as soon as we
show the existence of $N$), we need to prove the weak
   type of the maximal operator (i.e., (iii)). For a set
   $\Omega\subset\mathbb{R}^{n+1}_{+}$ and a function $u$ defined in
   $\mathbb{R}^{n+1}_{+}$ we define the maximal operator
   $\mathcal{M}_{\Omega}u(x)=\sup_{\Omega_{x}}|u|$.  Hence,
$\mathcal{M}_{\Omega}(P_t*f)(x)=M_{\Omega}f(x)$.
   
   We  can, without loss of generality, assume that the function $F$
   is positive. 
   First we split  the kernel $K_{t}(x)=P_t*K(x)$ into two parts, the local
   part of the kernel and the tail:  
   \[
   K_{t}(x)=\left(\chi _{|x| <3\gamma (t)}
   +\chi _{|x| >3\gamma (t)}\right)K_{t}(x)
   =K_{1,t}+K_{2,t}.
   \]
   
   First we consider the tail, $K_{2,t}$. 
   We need to estimate the following
   \[
  ( K_{2,t}*F)(x+x'), \ \mathrm{where}\ 
   (x',t)\in\Omega\subset \{(y,t):\ |y|\leq\gamma(t)\}.
   \]
   Assuming, $|x'|\leq \gamma (t)$, we have 
   \begin{eqnarray*}
     (K_{2,t}*F)(x+x')&=&\int_{\{|y|>3\gamma (t)\}}K_{t}(y)F(x+x'-y)dy \\
     &=&\int_{\{|y+x'|>3\gamma (t)\}}K_{t}(y+x')F(x-y)dy\\
     &\leq&\int_{\mathbb{R}}
     K_{t}(y/2)F(x-y)dy.     
   \end{eqnarray*}
   We know that since $K$ is radially decreasing,
   the same is true for $K_{t}$, and the boundedness of 
   $\mathcal{M}_{\Omega}(K_{2,t}*F)$ then follows from Lemma
   2.2  in \cite{NRS82}. 
   
   We now turn to the local  part of the kernel; i.e.,  $K_{1,t}$. 
   Let $\omega_{k}$ be the part
   of $\omega$ whose points  have the second coordinate equal to
   $t_{k}$:  
   \[  
   \omega_{k}=\{x;(x,t_{k})\in \omega \}.
   \]
   Let 
   \[
   \Omega _{k}=(\omega _{k}+{\Gamma})\cap \{(x,t); x\in
   {\mathbb{R}},\ t_{k}\leq t\leq t_{k-1}\}
   \]
   for  $k>1$, and for  $k=1$  let
   $\Omega _{1}=\omega_{1}+\Gamma$.  
   Then $\Omega \subset \Gamma_{3t}\cup (\cup \Omega_{k})$. 
   We split the operator as
   \begin{eqnarray*}
     \label{tk}
     \mathcal{M}_{\Omega }(K_{1,t}*F)(x)&\leq& \sup _{k} \mathcal{M}_{\Omega
       _{k}}(K_{1,t}*F)(x)+\mathcal{M}_{\Gamma_{3t}}(K_{1,t}*F)(x)  \\ &=&\sup
     _{k}T_{k}F(x)+ \mathcal{M}_{\Gamma_{3t}}(K_{1,t}*F)(x),
   \end{eqnarray*}
   where $T_{k}F(x)=\mathcal{M}_{\Omega
       _{k}}(K_{1,t}*F)(x)$.

   To use \lref{kbit} we need uniform weak type (1,1) estimates for
   the operators $T_{k}$, and they also have to fit the terminology of
   \lref{kbit}, which we will do below. 
   The main advantage of the lemma is that we can assume $t$ is in a
   fixed interval, away from $0$.

   To obtain the weak type (1,1) estimate, we first consider the part
   of $\Omega _{k}$ which 
   lies between the levels  $t_{k}$ and $2r(t_{k})$, which is
   $\Omega _{k}^{1}=\{(x,t)\in \Omega _{k};t_{k}<t<2r(t_{k})\}$.
   This part, $\Omega _{k}^{1}$, consists of $N(t_{k})$ non-tangential
   cones with vertices at the points 
   $(x_{i}^k,t_{k}),\ i=0,\dots ,N(t_{k})-1$. Let 
   $\Omega^{1}_{k,i}=((x_{i}^k,t_{k})+\Gamma)\cap \Omega_{k}^{1}$, 
   for
   $i=0,\dots ,N(t_{k})$, where we define $x^k_{N(t_{k})}=0$.
   Then, 
   \begin{eqnarray}  \label{kb8}    
        \|\mathcal{M}_{\Omega _{k}^{1}}(K_{1,t}*F)\|_{1,\infty }
      & \leq& \| \sup _{0\leq
         i\leq N(t_{k})-1} \ \ \mathcal{M}_{\Omega ^{1}_{k,i}}
       (K_{1,t}*F)\|_{1,\infty } \nonumber \\
       &\leq &\sum _{i=0}^{N(t_{k})-1}
       \| \mathcal{M}_{\Omega ^{1}_{k,i}}( K_{1,t}*F)
       \|_{1,\infty }\\
& \leq& N(t_{k})
       \|\mathcal{M}_{\Omega ^{1}_{k,N(t_{k})}}(K_{1,t}*F)
       \|_{1,\infty }\nonumber.
    \end{eqnarray}
   The last inequality follows from translation invariance. 
   The operator needs to be bounded uniformly in $k$ so we 
   need to see that the factor $N(t_{k})$ does not cause any problem.
   To proceed, we make a dyadic decomposition of the kernel
   $K_{1,t}$, and we get ($F$ is positive),
   \begin{eqnarray*}
     (K_{1,t}(x)
     \chi _{|x|<3\gamma (t)})*F&\leq & 
     \sum _{k=1}^{[C\log \gamma (t)/t]}(K_{1,t}(2^{k-1}t) \chi
     _{|x| <2^{k}t})*F
     \\ 
     &\leq  &
     C\sum _{k=1}^{[C\log \gamma
       (t)/t]}(2^{k-1}t)(K_{1,t}(2^{k-1}t) 
     \frac{1}{2^{k}t}\chi   _{|x| <2^{k}t})*F\\
     &\leq&
     C\sum _{k=1}^{[C\log \gamma
       (t)/t]}(2^{k-1}t)K_{1,t}(2^{k-1}t) 
     MF(x)\\
&\leq&
     CMF(x)\int_{t}^{\gamma (t)}K_{1,t}(x)dx,
   \end{eqnarray*}
   where $MF(x)$ is the usual Hardy-Littlewood maximal function.
   In order to bound \eref{kb8} uniformly in $k$, we
   must find a  bound on the integral times $N(t_{k})$, 
   for which we replace the limits with the smallest (respectively the largest) $t$
   allowed; i.e.,
   \[
   N(t_{k})\int_{t_{k}}^{\gamma (t_{k-1})}K_{1,t}(x)dx \leq
   N(t_{k})\|P_{t}\|_{L^{1}}\int_{t_{k}}^{\gamma (t_{k-1})}K(x)dx.
   \]
   The remaining integral in the right hand side decreases to $0$ as
   $t_{k}$ tends to $0$, i.e. $k\to\infty$, since both limits in the
   integral then tend to $0$ as $k\to\infty$. 
   Thus by choosing $N(t)$ to increase slowly
   enough to $\infty$ as $t\to 0$, we can bound the above
   expression uniformly in $k$.  This is the final restriction on
   $N(t)$. We have thus seen that it is always possible to find an
   unbounded $N(t)$ satisfying the above restrictions, as long as the
   kernel $K$ satisfies our assumptions. 

   Thus we can estimate the maximal operator by the usual 
   Hardy-Littlewood maximal function, which gives the weak type (1,1) for
   the operator $F\mapsto \mathcal{M}_{\Omega _{k}^{1}}(K_{1,t}*F)$
   uniformly in $k$, 
   if $k>1$.
   
   For the rest of $\Omega _{k}$, i.e.  if $2r(t_{k})\leq t \leq
   t_{k-1}$, then  $\Omega _{k}(t)$ consists of one interval, and we know
   from above that for this region we have the correct bound on the size. 

   Hence,  the  weak type (1,1) of the operator
   \begin{eqnarray*}     \label{kv}
     F\mapsto 
     \sup_{\substack{t>2r(t_{k}) \\ (x,t)\in \Omega_{k} }}
    ( K_{1,t}*F)(x),
   \end{eqnarray*}
   follows, since  the region we take the supremum  over is contained in a
   tangential region, of the right sort. This completes the proof of
   the uniform 
   weak type (1,1) of $T_{k},\ k>1$.

   The weak type (1,1) for $\mathcal{M}_{\Omega _{1}}(K_{1,t}*F)$
   follows by the same methods. 
   First take that part of $\Omega _{1}$ which lies between the levels
   $t_{1}$ and $2r(t_{1})$. Again, we will get a similar expression as
   above and this can be dealt with the same way.
   When $t>2r(t_{1})$, the region
   $\Omega _{1}$ is contained in  the tangential region
   $\Gamma_{K}$, and the weak   type (1,1) is proved. 

   Finally, we must check that our operators can be defined 
   as in \lref{kbit}, and that they satisfy the assumptions of the lemma.
   Let the index set $I_{k}$ be equal to $\Omega _{k}$ and set for
   $s=(x',t)\in I_{k}$,    
   \[
   K_{s}(x)=K_{1,t}(x+x'). 
   \]
   Then $T_{k}F(x)=\sup _{s\in I_{k}}(K_{s}*F)(x)$.
   To estimate the support of $K_{s}=K_{(x,t)}$, we see that the support
   is largest when $t=t_{k-1}$, which is the largest $t$ in the index set
   $I_{k}$. The support of $K_{1,t_{k-1}}$ is contained in the set
   $\{x;|x|\leq \gamma(t_{k-1})\}$, hence we can   bound  
   the support of $K_{s},\ s\in I_{k},$ taking $\gamma _{k}=3\gamma(t_{k-1})$.   
   If we take $N=2$,  then  we can bound the integral of $K_{s}^{*}$
   uniformly in $s\in\cup I_{k}$.
   With an $x$
   outside the support of $K_{s}$, we need only increase the support of
   the kernel $K_{1,t}$.   If  $s\in I_{k}$, using \eref{con} we obtain:
   \begin{eqnarray*}  
     \int_{0}^{\infty} K_{s}^{*}(x)dx&\leq&
     \int_{0}^{\gamma_{k+2}} K_{s}^{*}(x)dx+
     \int_{\gamma_{k+2}}^{\infty} K_{s}^{*}(x)dx
      \\& \leq&
     \int_{0}^{3\gamma (t_{k+1})} K_{t}(0)dx+
     \int_{3\gamma (t_{k+1})}^{\infty} K_{t}(x-3\gamma (t_{k+1}))dx
      \\ &\leq&
     3\frac{\gamma (t_{k+1})}{r(t_{k})}+     
     \int_{0}^{\infty} K_{t}(x)dx\leq
     3\frac{t_{k}}{r(t_{k})}+     
     \|K_{t}\|_{L^{1}},
   \end{eqnarray*}
   and from \eref{tang} it follows that   
   this expression is uniformly bounded  in $k$ for all 
   $s\in \cup I_{k}$.
   \lref{kbit} now gives the weak
   type (1,1) for $\sup _{k}T_{k}$, and hence for $\mathcal{M}_{\Omega }K_{1,t}$. 
   
   Finally, we have proved a weak type estimate for  
   both $\mathcal{M}_{\Omega }K_{1,t}$ and $\mathcal{M}_{\Omega
   }K_{2,t}$, and we have finished the proof of the theorem.  
\end{proof}

\bigbreak

\end{document}